# The distribution of the primes

*Cristiano Husu*

*Department of Mathematics*

*University of Connecticut*

**(A note to the readers: due to the many applications of the primes, the author trusts that the readers will handle the content with the requisite responsibility.)**

*To my father and mother*

## Abstract

This paper contains the extension of Euler's pentagonal number theorem to the sequence of the number of integer partitions with all parts equal and a visualization of the algorithm embodied in the theorem. There is a stunning identity at the end of the paper that implies that *the distribution of the primes is just a specific detail of the application of the sequence of the number of partitions (applied as Euler's inverse matrix) on a specific sequence, the $\sigma$-sequence, that we introduce and easily compute in the present paper*. Both the structure and the distribution of primes are determined just by the pentagonal number theorem sequence (Euler's sequence $\{1, -1, -1, 0, 0, 1, 0, 1, 0, \ldots\}$) and the pentagonal number themselves. Another striking feature of the algorithm, here, is that it detects new primes as quickly as Euler's algorithm computes the partition function. As for a quick summary of the content of the paper, the generating function of the sequence of the numbers of integer partitions with all parts equal, as a *series of rational functions*, is reassembled into an identity structured as Euler's pentagonal number theorem. The identity is composed of the series of rational functions itself, Euler's function $(1-x)(1-x^2)(1-x^3)\ldots$, and the $\sigma$-function, the result of the product of the series and Euler's function. Other than Euler's pentagonal number theorem, the other main idea, here, is to visualize the infinite *correction factors* contained in the series of rational functions as a matrix, the $\sigma$-matrix and the resulting $\sigma$-sequence (the sequence of the sums of the raising diagonals of the matrix) that takes the place of the sequence $\{1, 0, 0, 0, \ldots\}$ in Euler's algorithm. The $\sigma$-matrix *weaves in several different directions (vertically, backward horizontally, along the main diagonal and the other diagonal)* infinite copies of Euler's sequence. The $\sigma$-matrix is divided in three triangular areas, the upper triangular portion and the two triangular halves (connected to the origin) of the lower triangular portion. It is only in the lowest triangular portion of the $\sigma$-matrix that the correction factors of Euler's function create cluster of integers $k$ with the same sign and/or with $|k| \geq 2$, and these clusters, in turn, together with Euler's sequence, create a rich pentagonal tapestry of the distribution of the number of factors of positive integers.

## INTRODUCTION

The present work is a very short break in several years of study of various kinds of partitions with distinct parts and is loosely based on the author's previous work [H3] and [H4], which contains several applications of Euler's pentagonal number theorem and of another specialization of the Jacobi triple product. Conceptually, this paper is based on the creative way Euler found to manipulate functions and partition identities. It is also based on the work in [A1], [A2],[DL1], [DL2], [H1], [H2], [LW1], [LW2] and [TX]), for it is through set of books and papers that the author learned to evaluate and visualize various multiplications and divisions of generating functions by *correction factor*s (factors of Euler's function $(1-x)(1-x^2)(1-x^3)\ldots$). We refer the reader to the introduction of [H3] or [H4] for further background information.

The main object of the present work is *the generating function of the sequence of the numbers of integer partitions with all parts equal.* (The author does not like to use the words "number of divisors", as such words single out specific partitions and hide some of the correlations with other partitions.). This generating function is tautologically set to equal itself, as a *series of the rational functions* that filter the partition by the number of parts (or, equivalently, by the size of the parts). The resulting identity is multiplied on both sides by Euler's function, and reveals the *pentagonal number theorem for the number of partitions of* n*, denoted* $\rho(n)$*, with all parts equal* as the identity

$$\left[\prod_{m=1}^{\infty}(1-x^m)\right]\left[\sum_{n=0}^{\infty}\rho(n)x^n\right]=\sum_{h=0}^{\infty}x^h\prod_{m=1,m\neq h}^{\infty}(1-x^m),$$

(using the formal variable $x$ in place of the usual $q$). We name the right side of the above identity the $\sigma$-function.

The rest of the paper is a study of the $\sigma$-function along the structure of the pentagonal number theorem algorithm, the sole differences being that our algorithm, here, computes our desired generating function in place of Euler's computation of the partition function, and our algorithm requires the power series of the $\sigma$-function, whereas Euler's algorithm just has the constant function 1 as the product of Euler's function and the partition function.

Our computation of the power series of the $\sigma$-function, is visualized with a matrix that *weaves* infinite copies of the pentagonal number theorem sequence (the $\sigma$-matrix).

As for a detailed description of the content of the three sections of the paper, he first section contains the notation that we use and our version of the pentagonal number theorem. The second section is still written in the way the author uncovered the structure of the $\sigma$-matrix, the matrix visualization of the result of the multiplication of the Euler's function with the correction factors

$$\left(1+x^j+x^{2j}+x^{3j}+\cdots\right).$$

Finally, the third section contains the pentagonal algorithm, a visual comparison of Euler's matrix (the matrix of Euler's algorithm) and the $\sigma$-matrix. At the end of the paper, we show how the author, in search for an explicit formula for the primes, found out that the *$\sigma$-sequence* (the sum of the raising diagonals of the $\sigma$-matrix) is the missing link between the partition function and the main object of this work (identity (3.3), below).

We conclude the introduction with a short philosophical note. The fact that the distribution of the primes is as a weaving of the pentagonal number theorem sequence should not come as a surprise, because the multiplication of integers is just a useful place holder for the adding of partitions with all terms equal. Nature does not add with one hand and multiply with another, it just adds things up. It adds up when it puts things together and when it takes them apart.

A slightly different version of this paper is available at https://arxiv.org/pdf/xxxx.xxxxx. The author is very thankful to his wife, K. M. Johnson, an artist without a background in mathematics, for her guide and support during the frantic few days when we were, in a rough way we phrased it, "the only two people that know that the primes have only two ingredients: the outputs (1, 0 and −1) of the pentagonal number theorem sequence and the pentagonal numbers themselves". Without her colorful drawing of the first 102 rows and 102 columns of the $\sigma$-matrix, it would have taken much longer to write down the second section of this work.

# 1 THE GENERATING FUNCTION OF SEQUENCE OF THE NUMBER OF INTEGER PARTITIONS WITH ALL PARTS EQUAL AND THE $\sigma$-FUNCTION

For each integer $n$, let $\{\rho(n)\}$ be the sequence of the number of partitions of $n$ with all parts equal: set $\rho(n) = 0$, for $n < 0$, $\rho(0) = 1$, and, for each positive integer $n$, set $\rho(n)$ to be the number of partitions with parts $p_k$ with $k \geq 1$ and

$$p_1 = p_2 = p_3 = \cdots = p_k \geq 1. \tag{1.1}$$

Throughout the present work, we refer to the power series in Euler's pentagonal number theorem, the series $\sum_{i=-\infty}^{\infty}(-1)^i x^{\frac{3i^2}{2}+\frac{i}{2}} = 1 + \sum_{i=1}^{\infty}(-1)^i (x^{\frac{3i^2}{2}+\frac{i}{2}} + x^{\frac{3i^2}{2}-\frac{i}{2}})$, as the *pentagonal number theorem series*, and we refer to the underlying sequence, here (to be consistent with the related work and notation) denoted $\{\sigma_{n0}\}$ and $\{\sigma_0(n)\}$,

$$\{\sigma_{n0}\} = \{\sigma_0(n)\} = \{1, -1, -1, 0, 0, 1, 0, 1, 0, 0, 0, 0, -1, 0, 0, -1, \dots\}$$

for $n = 0, 1, 2, 3, \dots$, as the *pentagonal number theorem sequence.* To avoid any computational confusion, we also set $\sigma_{n0} = \sigma_0(n) = 0$, for $n < 0$.

As in the previous work [H3] and [H4], we use the formal variable $x$ in place of the customary $q$ to write down generating functions as $q$-functions. For the sequence $\{\rho(n)\}$, we write the generating function of the number of partitions of $n$ with all parts equal as

$$\sum_{n=0}^{\infty} \rho(n) x^n. \tag{1.2}$$

Note that both the filtration with respect to the number of parts $k$ in each partition, and the other filtration with respect to the size of the parts $|p_k|$ of the sequence $\{\rho(n)\}$ can be written as the same series of rational functions

$$\sum_{n=0}^{\infty} \rho(n)x^n = 1 + \frac{x}{1-x} + \frac{x^2}{1-x^2} + \frac{x^3}{1-x^3} + \cdots = 1 + \sum_{k=1}^{\infty} \frac{x^k}{1-x^k}, \tag{1.3}$$

as, for each $n \geq 0$, the partitions that add up to $\rho(n)$ can be analyzed in different equivalent ways. Multiplying (1.3) by Euler's function $\prod_{m=1}^{\infty} (1 - x^{\,m})$ yields the identity

$$\left[\prod_{m=1}^{\infty} (1-x^m)\right]\left[\sum_{n=0}^{\infty} \rho(n)x^n\right] = \sum_{h=0}^{\infty} x^h \prod_{m=1, m\neq h}^{\infty} (1-x^m). \tag{1.4}$$

We will refer to the function on the right side of (1.4) the *σ-function*.

The identity (1.4) has the same structure of Euler's partition function algorithm (cf.

$$\left[\prod_{m=1}^{\infty} (1-x^m)\right]\left[\sum_{n=0}^{\infty} p(n)x^n\right] = 1$$

with the partition function in place of (1.2) and the constant 1 in place of the $\sigma$-function). Comparing the two above identities reveals the rest of our work, here: (a) the computation of the *σ-function,* and (b) the pentagonal algorithm to compute $\rho(n)$ in terms of the sequence of the pentagonal number theorem. The next two sections correspond, respectively, to these two parts.

## 2 THE $\sigma$-FUNCTION AND THE $\sigma$-MATRIX

**Definition 2.1.** For each integer $j \geq 0$, we define *j-laced σ-sequences* $\{\sigma_j(n)\}$ as follows. Set $\sigma_j(n) = 0$, for $n < j$, and, on the other hand, for $0 \leq j \leq n$, set $\sigma_j(n)$ to be the coefficient of $n$ in $\prod_{m=1,m\neq j}^{\infty} (1 - x^m)$. □

We use the term *j-laced* as a visualization of the structure of the product of the series $1 + x^j + x^{2j} + x^{3j} + \cdots$ *(correction factors)* and Euler's function:

$$\prod_{m=1, m\neq j}^{\infty} (1-x^m) = \frac{1}{1-x^j}\prod_{m=1}^{\infty} (1-x^m) = (1 + x^j + x^{2j} + x^{3j} + \cdots)\prod_{m=1}^{\infty} (1-x^m) =$$

$$= \left(1 + x^j + x^{2j} + x^{3j} + \cdots\right)\sum_{i=-\infty}^{\infty} (-1)^i x^{\frac{3i^2}{2}+\frac{i}{2}} = \sum_{i=-\infty}^{\infty} (-1)^i x^{\frac{3i^2}{2}+\frac{i}{2}} + \sum_{i=-\infty}^{\infty} (-1)^i x^{\frac{3i^2}{2}+\frac{i}{2}+j} + \tag{2.1}$$

$$+ \sum_{i=-\infty}^{\infty} (-1)^i x^{\frac{3i^2}{2}+\frac{i}{2}+2j} + \sum_{i=-\infty}^{\infty} (-1)^i x^{\frac{3i^2}{2}+\frac{i}{2}+3j} + \dots .$$

□

(Note that the identity (2.1) is the pentagonal number theorem sequence counterpart of the familiar sifting of the multiplicative structure of integers.)

Lacing seems an appropriate word as, for any $n \geq 0$, the coefficient of $x^n$ in (1.5) is the sum of the coefficients of $x^n$, $x^{n-j}$, $x^{n-2j}$, $x^{n-3j}$, ... in the pentagonal number theorem sequence theorem (coefficients that are, therefore, "laced up" by the sum). Note that all these sums have a finite number of nonzero terms.

**Definition 2.2.** Define the *σ-sequence* $\{\sigma(n)\}$, for $n \geq 0$, to be the sequence of the coefficients of $n$ in the right side of the identity (1.4), that is the coefficients of $n$ in the expression

$$\prod_{m=1}^{\infty}(1-x^m)+x\prod_{m=1,m\neq 1}^{\infty}(1-x^m)+x^2\prod_{m=1,m\neq 2}^{\infty}(1-x^m)+x^3\prod_{m=1,m\neq 3}^{\infty}(1-x^m)+\cdots.$$

□

Note that the $\sigma$-sequence is the sequence underlying the power series expansion of the

$$\sum_{h=0}^{\infty}x^h\prod_{m=1,m\neq h}^{\infty}(1-x^m)=\sum_{n=0}^{\infty}\sigma(n)x^n,$$

and is also the sum of the $j$-laced $\sigma$-sequences as

$$\sigma(n)=\sum_{j=0}^{n}\sigma_j(n-j).$$

**Definition 2.3.** Define the *σ-matrix* be the infinite matrix $\{\sigma_{ij}\}$, $i,j=0,1,3,\dots$, in the following way. Set $\sigma_{ij}=0$, when one, or both, the sub-indices are negative, $i<0$ or $j<0$. The matrix is nonzero only in the fourth quadrant. The first nonzero column, for $j=0$, is the pentagonal number theorem sequence $\sigma_{i0}=\sigma_0(i)$ the coefficient of $x^i$ in $\prod_{m=1}^{\infty}(1-x^m)$ (at the beginning of the previous section), the first row is just a row of ones $\sigma_{0j}=\sigma_j(0)=1$, and the rest of the matrix, for $i,j>0$, is *laced* in the triangular recursive way

$$\{\sigma_{ij}\}=\{\sigma_{i0}\}+\{\sigma_{i-j,j}\}. \tag{2.2}$$

□

Note that $\sigma_{ij}=\sigma_j(i)$, for all $i,j\geq 0$, the correction factor $1+x^j+x^{2j}+x^{3j}+\cdots$ in (2.1) implies, by induction, that $\sigma_j(i)=\sigma_0(i)+\sigma_0(i-j)+\sigma_0(i-2j)+\sigma_0(i-3j)+\cdots=\sigma_0(i)+\sigma_j(i-j)$ as $\sigma_j(i-j)=\sigma_0(i-2j)+\sigma_0(i-3j)+\cdots$, and so forth.

In the extended version of the present paper, [H5], the reader can see the first few terms (the terms in the first 102 by 102 square for $n=101$) of the $\sigma$-matrix, computed with pencil and paper. Several properties of both $\sigma$-matrix and the $\sigma$-function are easy to see and prove. Some of the following properties make it easier to compute both the $\sigma$-matrix and the $\sigma$-function for larger numbers of rows or columns.

**Property 2.4.** $\sigma_{ij} = \sigma_{i0}$, if $i > j$, so that the terms in the upper triangular matrix (in the fourth quadrant) are just copies of the corresponding terms (the projections on the vertical axis) of the pentagonal number sequence.

**Property 2.5.** $\sigma_{ij} \geq 0$, if $i = j$.

**Property 2.6.** $|\sigma_{ij}| \geq 2$ implies that $i$ is the sum of a (generalized) pentagonal number and a multiple of $j$.

The next property is the main clue that reveals the structure of the $\sigma$-matrix.

**Property 2.7.** $\sigma_{i,j-2} = -1$, $\sigma_{i,j-1} = -1$, $\sigma_{ij} = 1$, if $i \neq 3$ and $i$ is not multiple of a (generalized) pentagonal number.

To be more specific with the details of the pentagonal algorithm, here and in the next section, we will use the following descriptions of pentagonal numbers.

**Definition 2.8. (a)** Denoting the (generalized) pentagonal numbers as $P_n = \frac{3n^2}{2} - \frac{n}{2}$, with $n = 0, \pm1, \pm2, \pm3, \ldots$ , we can set the *even-indexed (EI) pentagonal numbers* to be the numbers $P_{2k} = 6k^2 - k$, with $k = 0, \pm1, \pm2, \pm3, \ldots$ , and the *odd-indexed (OI) pentagonal numbers* to be $P_{2k+1} = 6k^2 + 5k + 1$, again with $k = 0, \pm1, \pm2, \pm3, \ldots$, so that the counterimages of 1, and respectively, $-1$, are the even-indexed pentagonal numbers, and, respectively, the odd-indexed pentagonal numbers.

**(b)** We also define the $EI-$ (respectively, $EI+$) *pentagonal numbers, the even-indexed negative (respectively, positive) pentagonal numbers,* to be the numbers of the form $P_{2k}^{-} = 6k^2 - k$ (respectively, $P_{2k}^{+} = 6k^2 + k$), for $k = 0, 1, 2, 3, \ldots$ , and the $0I-$ (respectively, $OI+$) *pentagonal numbers, the odd-indexed negative (respectively, positive) pentagonal numbers,* to be $P_{2k+1}^{-} = 6k^2 + 5k + 1$, (respectively, $P_{2k+1}^{+} = 6k^2 + 7k + 2$), again for $k = 0, 1, 2, 3, \ldots$ . □

The following property is a refinement of Property 2.6 and provides a deep insight of the way the pentagonal number theorem sequence shapes the distribution of the partitions with an equal number of parts, as the larger values and, respectively, the smaller values of $\sigma_{ij}$ *lump up* around the rows corresponding to the even-indexed pentagonal numbers and, respectively, the odd-indexed pentagonal numbers, and these lumps will determine large changes in the number of such partitions.

**Property 2.9.** $\sigma_{ij} \geq 2$ implies that $i$ is the sum of an even-indexed pentagonal numbers and a multiple of $j$. On the other hand, $\sigma_{ij} \leq -2$ implies that $i$ is the sum of an even-indexed pentagonal numbers and a multiple of $j$.

**Property 2.10.** The terms of the $\sigma$-sequence $\{\sigma(n)\}$ are the sum of the terms of the diagonals, increasing at the 45° angle, in the $\sigma$-matrix starting at the top left corner, so, for example,

$\{\sigma(0)\} = \sigma_{00} = 1$,

$\{\sigma(1)\} = \sigma_{10} + \sigma_{01} = 1 - 1 = 0$,

$\{\sigma(2)\} = \sigma_{20} + \sigma_{11} + \sigma_{02} = 1 + 0 - 1 = 0$,

$\{\sigma(3)\} = \sigma_{30} + \sigma_{21} + \sigma_{12} + \sigma_{03} = 1 - 1 - 1 + 0 = -1$, and so forth.

**Definition 2.11.** Define the upper $\sigma$-sequence to be $\sigma_u(n) = \sum_{j=0,n>2j}^{n} \sigma_j(n-j)$, and the lower $\sigma$-sequence to be $\sigma_l(n) = \sum_{j=0,n\leq 2j}^{n} \sigma_j(n-j)$.

**Property 2.12.** The upper $\sigma$-sequence is the sequence of the partial sums of the pentagonal number theorem sequence, which means, using any fixed $k = 0, 1, 2, 3, \dots$,

$\sigma_u(n) = 1$, for $6k^2 + k \leq n < 6k^2 + 5k + 1$ (for all integers $n$ between any $EI+$ pentagonal number and the next pentagonal number),

$\sigma_u(n) = -1$, for $6k^2 + 7k + 2 \leq n < 6k^2 + 11k + 5$ (for all $n$ between any $OI+$ pentagonal number and the next pentagonal number), and

$\sigma_u(n) = 0$, for $6k^2 - k \leq n < 6k^2 + k$, and for $6k^2 + 5k + 1 \leq n < 6k^2 + 7k + 2$ (for all $n$ between any two consecutive pentagonal numbers with the same sign).

**Definition 2.13.** Define the *upper triangular $\sigma$-matrix* to be the infinite triangular matrix $\{v_{ij}\} = \{\sigma_{ij}\}$, $i, j = 0, 1, 3, \dots$, and $j > i$. Define the *λ-matrix,* the *lower triangular $\sigma$-matrix* to be the infinite matrix $\{\lambda_{ij}\} = \{\sigma_{ij}\}$, $i, j = 0, 1, 3, \dots$, and $j \leq i$. Note that the main diagonal of the $\sigma$-matrix is not included in the upper triangular matrix, but it is included in the $\Lambda$-matrix. We also define the *$\lambda_1$-matrix* to be the matrix $\{\lambda_{ij}\} = \{\sigma_{ij}\}$, $i, j = 0, 1, 3, \dots$, and $j \leq i \leq 2i$, and the *$\lambda_2$-matrix* to be the matrix $\{\lambda_{ij}\} = \{\sigma_{ij}\}$, $i, j = 0, 1, 3, \dots$, and $i > 2i$. Note that the line $i = 2i$, the *correction border,* is included the *$\lambda_1$-matrix,* but not in the $\lambda_2$-matrix. We do refer to the $\lambda_2$-matrix as the *correction zone* of the *$\sigma$-matrix* for it is in this area that the correction factors in (2.1) weave the pentagonal number theorem sequence. □

Once an $n+1$ times $n+1$ portion of the $\Lambda$-matrix is computed (Property 2.4 implies that the upper triangular $\sigma$-matrix does not need to be computed, it is just set in place for the triangulation (2.2)), we can quickly see the $\sigma$-sequence (half of which predetermined by properties 2.4, 2.10 and 2.12) and add up the power series of the $\sigma$-function (as in right hand side of (1.6)) as, for example, for $n = 101$,

$$\sum_{n=0}^{\infty} \sigma(n)x^n = 1 - x^3 - x^4 - 2x^5 - x^7 + x^8 + 2x^9 + x^{10} + 2x^{11} + 3x^{12} + x^{13} - x^{14} + 3x^{15} +$$

$$-2x^{16} - 2x^{17} - x^{18} - 3x^{19} - 2x^{20} - 2x^{21} - 5x^{22} - x^{24} + x^{25} - 3x^{26} + \cdots$$

and the underlying $\sigma$-sequence as

$$\{\sigma(n)\} = \{1, 0, 0, -1, -1, -2, 0, -1, 1, 2, 1, 2, 3, 1, -1, 3, -2, -2, -1, -3, -2, -2, -5, 0, -1, 1, -3, \dots\}$$

## 3 THE PENTAGONAL ALGORITHM AND THE STRUCTURE OF THE σ-MATRIX

As in Euler's Pentagonal Number Theorem ([A2], Theorem 1.6 and Corollaries 1.7, 1.8). We multiply the factors on the left side of (1.4) as

$$\left[\sum_{i=-\infty}^{\infty}(-1)^i x^{\frac{3i^2}{2}+\frac{i}{2}}\right]\left[\sum_{n=0}^{\infty}\rho(n)x^n\right] = \sum_{i=-\infty}^{\infty}\sum_{n=0}^{\infty}(-1)^i\rho(n)x^{n+\frac{3i^2}{2}+\frac{i}{2}} =$$

$$= \sum_{i=-\infty}^{\infty}\sum_{m=0}^{\infty}(-1)^i\rho\left(m-\frac{3i^2}{2}-\frac{i}{2}\right)x^m =$$

so that the identity (1.4) can be written as

$$\sum_{i=-\infty}^{\infty}\sum_{m=0}^{\infty}(-1)^i\rho\left(m-\frac{3i^2}{2}-\frac{i}{2}\right)x^m = \sum_{n=0}^{\infty}\sigma(n)x^n, \tag{3.1}$$

and the pentagonal algorithm as

$$\sum_{i=-\infty}^{\infty}(-1)^i\rho\left(n-\frac{3i^2}{2}-\frac{i}{2}\right)x^m = \sigma(n) \tag{3.2}$$

(a finite sum since $\rho(n) = 0$, for $n < 0$).

To see the structure of the $\sigma$-matrix, fix a nonnegative integer $n$, unfold (3.2) as

$$\begin{aligned}
&\rho(0) &&= \sigma(0)\\
&-\rho(0)+\rho(1) &&= \sigma(1)\\
&-\rho(0)-\rho(1)+\rho(2) &&= \sigma(2)\\
&\qquad -\rho(1)-\rho(2)+\rho(3) &&= \sigma(3)\\
&\qquad\qquad -\rho(2)-\rho(3)+\rho(4) &&= \sigma(4)\\
&\rho(0)+ \qquad\qquad -\rho(3)-\rho(4)+\rho(5) &&= \sigma(5)\\
&\qquad \rho(1)+ \qquad\qquad -\rho(4)-\rho(5)+\rho(6) &&= \sigma(6)\\
&\rho(0)+ \qquad \rho(2)+ \qquad\qquad -\rho(5)-\rho(6)+\rho(7) &&= \sigma(7)\\
&\qquad \rho(1)+ \qquad \rho(3)+ \qquad\qquad -\rho(6)-\rho(7)+\rho(8) &&= \sigma(8)
\end{aligned}$$

$$\rho(2) + \qquad \rho(4) + \qquad -\rho(7) - \rho(8) + \rho(9) \qquad = \sigma(9)$$

$$\ldots\ldots\ldots + \rho(n-7) \qquad + \rho(n-5) \qquad - \rho(n-2) - \rho(n-1) + \rho(n) = \sigma(n),$$

write down the $(n+1) \times (n+1)$ matrix of the above pentagonal algorithm, the Euler matrix, since it is the same for both algorithms. For example, for $n = 25$, the matrix fits, here, as

| | | | | | | | | | | | | | | | | | | | | | | | | |
|---|---|---|---|---|---|---|---|---|---|---|---|---|---|---|---|---|---|---|---|---|---|---|---|---|
| 1 | 0 | 0 | 0 | 0 | 0 | 0 | 0 | 0 | 0 | 0 | 0 | 0 | 0 | 0 | 0 | 0 | 0 | 0 | 0 | 0 | 0 | 0 | 0 | 0 |
| -1 | 1 | 0 | 0 | 0 | 0 | 0 | 0 | 0 | 0 | 0 | 0 | 0 | 0 | 0 | 0 | 0 | 0 | 0 | 0 | 0 | 0 | 0 | 0 | 0 |
| -1 | -1 | 1 | 0 | 0 | 0 | 0 | 0 | 0 | 0 | 0 | 0 | 0 | 0 | 0 | 0 | 0 | 0 | 0 | 0 | 0 | 0 | 0 | 0 | 0 |
| 0 | -1 | -1 | 1 | 0 | 0 | 0 | 0 | 0 | 0 | 0 | 0 | 0 | 0 | 0 | 0 | 0 | 0 | 0 | 0 | 0 | 0 | 0 | 0 | 0 |
| 0 | 0 | -1 | -1 | 1 | 0 | 0 | 0 | 0 | 0 | 0 | 0 | 0 | 0 | 0 | 0 | 0 | 0 | 0 | 0 | 0 | 0 | 0 | 0 | 0 |
| 1 | 0 | 0 | -1 | -1 | 1 | 0 | 0 | 0 | 0 | 0 | 0 | 0 | 0 | 0 | 0 | 0 | 0 | 0 | 0 | 0 | 0 | 0 | 0 | 0 |
| 0 | 1 | 0 | 0 | -1 | -1 | 1 | 0 | 0 | 0 | 0 | 0 | 0 | 0 | 0 | 0 | 0 | 0 | 0 | 0 | 0 | 0 | 0 | 0 | 0 |
| 1 | 0 | 1 | 0 | 0 | -1 | -1 | 1 | 0 | 0 | 0 | 0 | 0 | 0 | 0 | 0 | 0 | 0 | 0 | 0 | 0 | 0 | 0 | 0 | 0 |
| 0 | 1 | 0 | 1 | 0 | 0 | -1 | -1 | 1 | 0 | 0 | 0 | 0 | 0 | 0 | 0 | 0 | 0 | 0 | 0 | 0 | 0 | 0 | 0 | 0 |
| 0 | 0 | 1 | 0 | 1 | 0 | 0 | -1 | -1 | 1 | 0 | 0 | 0 | 0 | 0 | 0 | 0 | 0 | 0 | 0 | 0 | 0 | 0 | 0 | 0 |
| 0 | 0 | 0 | 1 | 0 | 1 | 0 | 0 | -1 | -1 | 1 | 0 | 0 | 0 | 0 | 0 | 0 | 0 | 0 | 0 | 0 | 0 | 0 | 0 | 0 |
| 0 | 0 | 0 | 0 | 1 | 0 | 1 | 0 | 0 | -1 | -1 | 1 | 0 | 0 | 0 | 0 | 0 | 0 | 0 | 0 | 0 | 0 | 0 | 0 | 0 |
| -1 | 0 | 0 | 0 | 0 | 1 | 0 | 1 | 0 | 0 | -1 | -1 | 1 | 0 | 0 | 0 | 0 | 0 | 0 | 0 | 0 | 0 | 0 | 0 | 0 |
| 0 | -1 | 0 | 0 | 0 | 0 | 1 | 0 | 1 | 0 | 0 | -1 | -1 | 1 | 0 | 0 | 0 | 0 | 0 | 0 | 0 | 0 | 0 | 0 | 0 |
| 0 | 0 | -1 | 0 | 0 | 0 | 0 | 1 | 0 | 1 | 0 | 0 | -1 | -1 | 1 | 0 | 0 | 0 | 0 | 0 | 0 | 0 | 0 | 0 | 0 |
| -1 | 0 | 0 | -1 | 0 | 0 | 0 | 0 | 1 | 0 | 1 | 0 | 0 | -1 | -1 | 1 | 0 | 0 | 0 | 0 | 0 | 0 | 0 | 0 | 0 |
| 0 | -1 | 0 | 0 | -1 | 0 | 0 | 0 | 0 | 1 | 0 | 1 | 0 | 0 | -1 | -1 | 1 | 0 | 0 | 0 | 0 | 0 | 0 | 0 | 0 |
| 0 | 0 | -1 | 0 | 0 | -1 | 0 | 0 | 0 | 0 | 1 | 0 | 1 | 0 | 0 | -1 | -1 | 1 | 0 | 0 | 0 | 0 | 0 | 0 | 0 |
| 0 | 0 | 0 | -1 | 0 | 0 | -1 | 0 | 0 | 0 | 0 | 1 | 0 | 1 | 0 | 0 | -1 | -1 | 1 | 0 | 0 | 0 | 0 | 0 | 0 |
| 0 | 0 | 0 | 0 | -1 | 0 | 0 | -1 | 0 | 0 | 0 | 0 | 1 | 0 | 1 | 0 | 0 | -1 | -1 | 1 | 0 | 0 | 0 | 0 | 0 |
| 0 | 0 | 0 | 0 | 0 | -1 | 0 | 0 | -1 | 0 | 0 | 0 | 0 | 1 | 0 | 1 | 0 | 0 | -1 | -1 | 1 | 0 | 0 | 0 | 0 |
| 0 | 0 | 0 | 0 | 0 | 0 | -1 | 0 | 0 | -1 | 0 | 0 | 0 | 0 | 1 | 0 | 1 | 0 | 0 | -1 | -1 | 1 | 0 | 0 | 0 |
| 1 | 0 | 0 | 0 | 0 | 0 | 0 | -1 | 0 | 0 | -1 | 0 | 0 | 0 | 0 | 1 | 0 | 1 | 0 | 0 | -1 | -1 | 1 | 0 | 0 |
| 0 | 1 | 0 | 0 | 0 | 0 | 0 | 0 | -1 | 0 | 0 | -1 | 0 | 0 | 0 | 0 | 1 | 0 | 1 | 0 | 0 | -1 | -1 | 1 | 0 |
| 0 | 0 | 1 | 0 | 0 | 0 | 0 | 0 | 0 | -1 | 0 | 0 | -1 | 0 | 0 | 0 | 0 | 1 | 0 | 1 | 0 | 0 | -1 | -1 | 1 |
| 0 | 0 | 0 | 1 | 0 | 0 | 0 | 0 | 0 | 0 | -1 | 0 | 0 | -1 | 0 | 0 | 0 | 0 | 1 | 0 | 1 | 0 | 0 | -1 | -1 |

(if the reader has a matrix C.A.S., try much larger numbers).

We place the above matrix next to the corresponding portion of the $\sigma$-matrix

| | | | | | | | | | | | | | | | | | | | | | | | | |
|---|---|---|---|---|---|---|---|---|---|---|---|---|---|---|---|---|---|---|---|---|---|---|---|---|
| 1 | 1 | 1 | 1 | 1 | 1 | 1 | 1 | 1 | 1 | 1 | 1 | 1 | 1 | 1 | 1 | 1 | 1 | 1 | 1 | 1 | 1 | 1 | 1 | 1 |
| -1 | 0 | -1 | -1 | -1 | -1 | -1 | -1 | -1 | -1 | -1 | -1 | -1 | -1 | -1 | -1 | -1 | -1 | -1 | -1 | -1 | -1 | -1 | -1 | -1 |
| -1 | -1 | 0 | -1 | -1 | -1 | -1 | -1 | -1 | -1 | -1 | -1 | -1 | -1 | -1 | -1 | -1 | -1 | -1 | -1 | -1 | -1 | -1 | -1 | -1 |
| 0 | -1 | -1 | 1 | 0 | 0 | 0 | 0 | 0 | 0 | 0 | 0 | 0 | 0 | 0 | 0 | 0 | 0 | 0 | 0 | 0 | 0 | 0 | 0 | 0 |
| 0 | -1 | 0 | -1 | 1 | 0 | 0 | 0 | 0 | 0 | 0 | 0 | 0 | 0 | 0 | 0 | 0 | 0 | 0 | 0 | 0 | 0 | 0 | 0 | 0 |
| 1 | 0 | 0 | 0 | 0 | 2 | 1 | 1 | 1 | 1 | 1 | 1 | 1 | 1 | 1 | 1 | 1 | 1 | 1 | 1 | 1 | 1 | 1 | 1 | 1 |
| 0 | 0 | 0 | 1 | -1 | -1 | 1 | 0 | 0 | 0 | 0 | 0 | 0 | 0 | 0 | 0 | 0 | 0 | 0 | 0 | 0 | 0 | 0 | 0 | 0 |
| 1 | 1 | 1 | 0 | 1 | 0 | 0 | 2 | 1 | 1 | 1 | 1 | 1 | 1 | 1 | 1 | 1 | 1 | 1 | 1 | 1 | 1 | 1 | 1 | 1 |
| 0 | 1 | 0 | 0 | 1 | 0 | -1 | -1 | 1 | 0 | 0 | 0 | 0 | 0 | 0 | 0 | 0 | 0 | 0 | 0 | 0 | 0 | 0 | 0 | 0 |
| 0 | 1 | 1 | 1 | 0 | 0 | 0 | -1 | -1 | 1 | 0 | 0 | 0 | 0 | 0 | 0 | 0 | 0 | 0 | 0 | 0 | 0 | 0 | 0 | 0 |
| 0 | 1 | 0 | 0 | -1 | 2 | 0 | 0 | -1 | -1 | 1 | 0 | 0 | 0 | 0 | 0 | 0 | 0 | 0 | 0 | 0 | 0 | 0 | 0 | 0 |
| 0 | 1 | 1 | 0 | 1 | -1 | 1 | 0 | 0 | -1 | -1 | 1 | 0 | 0 | 0 | 0 | 0 | 0 | 0 | 0 | 0 | 0 | 0 | 0 | 0 |
| -1 | 0 | -1 | 0 | 0 | -1 | 0 | 0 | -1 | -1 | -2 | -2 | -1 | -1 | -1 | -1 | -1 | -1 | -1 | -1 | -1 | -1 | -1 | -1 | -1 |
| 0 | 0 | 1 | 0 | 0 | 0 | 0 | 0 | 1 | 0 | 0 | -1 | -1 | 1 | 0 | 0 | 0 | 0 | 0 | 0 | 0 | 0 | 0 | 0 | 0 |
| 0 | 0 | -1 | 0 | -1 | 0 | -1 | 2 | 0 | 1 | 0 | 0 | -1 | -1 | 1 | 0 | 0 | 0 | 0 | 0 | 0 | 0 | 0 | 0 | 0 |
| -1 | -1 | 0 | -1 | 0 | 1 | -1 | -2 | 0 | -1 | 0 | -1 | -1 | -2 | -2 | 0 | -1 | -1 | -1 | -1 | -1 | -1 | -1 | -1 | -1 |
| 0 | -1 | -1 | 0 | 0 | -1 | 0 | -1 | 1 | 1 | 0 | 1 | 0 | 0 | -1 | -1 | 1 | 0 | 0 | 0 | 0 | 0 | 0 | 0 | 0 |
| 0 | -1 | 0 | 0 | 0 | -1 | 1 | 0 | -1 | 0 | 1 | 0 | 1 | 0 | 0 | -1 | -1 | 1 | 0 | 0 | 0 | 0 | 0 | 0 | 0 |
| 0 | -1 | -1 | -1 | -1 | 0 | 0 | 0 | -1 | 1 | 0 | 1 | 0 | 1 | 0 | 0 | -1 | -1 | 1 | 0 | 0 | 0 | 0 | 0 | 0 |
| 0 | -1 | 0 | 0 | 0 | 0 | 0 | 0 | 0 | -1 | 0 | 0 | 1 | 0 | 1 | 0 | 0 | -1 | -1 | 1 | 0 | 0 | 0 | 0 | 0 |
| 0 | -1 | -1 | 0 | 0 | 1 | -1 | 0 | -1 | -1 | 1 | 0 | 0 | 1 | 0 | 1 | 0 | 0 | -1 | -1 | 1 | 0 | 0 | 0 | 0 |
| 0 | -1 | 0 | -1 | 0 | -1 | -1 | 2 | 1 | -1 | -1 | 0 | 0 | 0 | 1 | 0 | 1 | 0 | 0 | -1 | -1 | 1 | 0 | 0 | 0 |
| 1 | 0 | 0 | 1 | 0 | 0 | 1 | -1 | 1 | 1 | -1 | 2 | 1 | 1 | 1 | 2 | 1 | 2 | 1 | 1 | 0 | 0 | 2 | 0 | 0 |
| 0 | 0 | 0 | 0 | 0 | 0 | 1 | -1 | 0 | 1 | 0 | -2 | 0 | 0 | 0 | 0 | 1 | 0 | 1 | 0 | 0 | -1 | -1 | 1 | 0 |
| 0 | 0 | 0 | -1 | 0 | 0 | 0 | 0 | 1 | -1 | 0 | -1 | 0 | 0 | 0 | 0 | 0 | 1 | 0 | 1 | 0 | 0 | -1 | -1 | 1 |
| 0 | 0 | 0 | 1 | 0 | 1 | 0 | 0 | -1 | 1 | 0 | 0 | -1 | 0 | 0 | 0 | 0 | 1 | 0 | 1 | 0 | 0 | -1 | -1 | 1 |

and **see** the structure of Euler's matrix in the $\sigma$-matrix:

**Remark 3.1.** In the above visualization of the $\sigma$-matrix, it is easy to see the structure of the $\sigma$-matrix, and how the three parts of the matrix affect the ups and downs of the $\sigma$-sequence. The upper triangular matrix $\{v_{ij}\}$ is just made of vertical copies of the pentagonal number sequence. The highlighted portion, the $\lambda_1$-matrix, is just made up of copies of horizontal copies of the pentagonal number sequence *written backward* (the yellow highlighted rows), starting on the main diagonal and ending on the (green) correction border, and *with one of the (infinite repeating) numbers on the right of the main diagonal, on that row, added to the sequence.* The portion on the left, the $\lambda_2$-matrix, is determined in the same way as in the second part but with additional action of the correction factors that, in that section of the paper add an additional layer to the matrix (cf. Property 2.6). □

We conclude our work, here, with the inverse matrix version of our algorithm.

Fix a nonnegative integer $n$, denote the first $n+1$ rows and the first $n+1$ columns of Euler matrix $E_{(n+1)\times(n+1)}$and denote $E^{-1}_{(n+1)\times(n+1)}$ its inverse matrix. Denote $\overrightarrow{\rho(n)}$, and, respectively, $\overrightarrow{\sigma(n)}$, the first $n+1$ terms of our desired sequence, and, respectively, the $\sigma$-sequence. Then, our algorithm can be written as

$$\overrightarrow{\rho(n)} = E^{-1}_{(n+1)\times(n+1)}\overrightarrow{\sigma(n)}. \tag{3.3}$$

Note that the matrix $E^{-1}_{(n+1)\times(n+1)}$is built just structured like the matrix $E_{(n+1)\times(n+1)}$, except for the partition function sequence $\{1,1,2,3,5,7,11,15,22,30,42,56,77,101,135,176,231,297,\ldots\}$, *written backward*, as

{… 297, 231, 176, 135, 101, 77, 56, 42, 30, 22, 15, 11, 7, 5, 3, 2, 1, 1}

*and showing up one term at a time in each new row* (as in the structure of $E_{(n+1)\times(n+1)}$) in place of Euler's sequence (for $n = 25$, in the highlighted whole lower triangular portion of Euler's matrix). In other words, the algorithm, written as in (3.3), is telling us that, for all nonnegative integers $n$, the structure of the sequence of the number of partitions with all equal parts is determined by the partition function, *structured as in the matrices* $E_{(n+1)\times(n+1)}$, applied to the $\sigma$-sequence. In particular, the $\sigma$-sequence embodies the relationship between the partition function and the distribution of the primes.